\documentclass[12pt]{amsart}
\usepackage{amssymb, amsmath}
\usepackage{verbatim} 
\usepackage{graphicx}
\usepackage{wrapfig}
\usepackage{caption}
\usepackage{subcaption}
\usepackage{color}
 \usepackage{array}
 \newcolumntype{o}{@{}>{{}}c<{{}}@{}}

\usepackage{colortbl}

\usepackage{enumitem}
 \usepackage{empheq}

\begin{document}
\newtheorem{thm}{Theorem}[section]
\newtheorem*{thm*}{Theorem}
\newtheorem{lem}[thm]{Lemma}
\newtheorem{prop}[thm]{Proposition}
\newtheorem{cor}[thm]{Corollary}
\newtheorem*{conj}{Conjecture}

\theoremstyle{definition}
\newtheorem*{defn}{Definition}
\newtheorem*{remark}{Remark}
\newtheorem{exercise}{Exercise}
\newtheorem*{exercise*}{Exercise}

 \numberwithin{equation}{section}

\newcommand{\rad}{\operatorname{rad}}

\newcommand{\Z}{{\mathbb Z}} 
\newcommand{\Q}{{\mathbb Q}}
\newcommand{\R}{{\mathbb R}}
\newcommand{\C}{{\mathbb C}}
\newcommand{\N}{{\mathbb N}}
\newcommand{\FF}{{\mathbb F}}
\newcommand{\fq}{\mathbb{F}_q}
\newcommand{\rmk}[1]{\footnote{{\bf Comment:} #1}}

\renewcommand{\mod}{\;\operatorname{mod}}
\newcommand{\ord}{\operatorname{ord}}
\newcommand{\TT}{\mathbb{T}}
\renewcommand{\i}{{\mathrm{i}}}
\renewcommand{\d}{{\mathrm{d}}}
\renewcommand{\^}{\widehat}
\newcommand{\HH}{\mathbb H}
\newcommand{\Vol}{\operatorname{vol}}
\newcommand{\area}{\operatorname{area}}
\newcommand{\tr}{\operatorname{tr}}
\newcommand{\norm}{\mathcal N} 
\newcommand{\intinf}{\int_{-\infty}^\infty}
\newcommand{\ave}[1]{\left\langle#1\right\rangle} 
\newcommand{\Var}{\operatorname{Var}}
\newcommand{\Prob}{\operatorname{Prob}}
\newcommand{\sym}{\operatorname{Sym}}
\newcommand{\disc}{\operatorname{disc}}
\newcommand{\CA}{{\mathcal C}_A}
\newcommand{\cond}{\operatorname{cond}} 
\newcommand{\lcm}{\operatorname{lcm}}
\newcommand{\Kl}{\operatorname{Kl}} 
\newcommand{\leg}[2]{\left( \frac{#1}{#2} \right)}  
\newcommand{\Li}{\operatorname{Li}}

\newcommand{\sumstar}{\sideset \and^{*} \to \sum}

\newcommand{\LL}{\mathcal L} 
\newcommand{\sumf}{\sum^\flat}
\newcommand{\Hgev}{\mathcal H_{2g+2,q}}
\newcommand{\USp}{\operatorname{USp}}
\newcommand{\conv}{*}
\newcommand{\dist} {\operatorname{dist}}
\newcommand{\CF}{c_0} 
\newcommand{\kerp}{\mathcal K}

\newcommand{\Cov}{\operatorname{cov}}
\newcommand{\Sym}{\operatorname{Sym}}

\newcommand{\Ht}{\operatorname{Ht}}

\newcommand{\E}{\operatorname{E}} 
\newcommand{\sign}{\operatorname{sign}} 
\newcommand{\meas}{\operatorname{meas}} 

\newcommand{\divid}{d} 

\newcommand{\GL}{\operatorname{GL}}
\newcommand{\SL}{\operatorname{SL}}
\newcommand{\re}{\operatorname{Re}}
\newcommand{\im}{\operatorname{Im}}
\newcommand{\res}{\operatorname{Res}}
\newcommand{\Ex}{\mathcal E} 

\title{Zeros of modular forms and Faber polynomials}
\author{Ze\'ev Rudnick}
\address{School of Mathematical Sciences, Tel Aviv University, Tel Aviv 69978, Israel}
\email{rudnick@tauex.tau.ac.il}
\date{\today}
\begin{abstract}
We study the zeros of cusp forms of large weight for the modular group, which have a very large order of vanishing at infinity, so that they have a fixed number $D$ of finite zeros in the fundamental domain. 
We show that for large weight the zeros of these forms cluster near $D$ vertical lines, with the zeros of a weight $k$ form lying at  height approximately $\log k$.  
This is in contrast to previously known cases, such as Eisenstein series,  where the zeros lie on the circular part of the boundary of the fundamental domain, or the case of cuspidal Hecke eigenforms where the zeros are uniformly distributed in the fundamental domain.

Our method uses the Faber polynomials.  We  show that for our class of cusp forms, the associated Faber polynomials, suitably renormalized, converge to the truncated exponential polynomial of degree $D$.

\end{abstract}
\maketitle
 
 \section{Introduction}
 
For an even integer $k\geq 0$  let $M_k$ be the space of modular forms of  weight $k$ for the full modular group $\SL(2,\Z)$.
Any $f\in M_k$ has an expansion
\[
f(\tau) = \sum_{n=0}^\infty a_f(n)q^n
\]
where   $q=e^{2\pi i \tau}$, $\tau\in \HH =\{ \tau:\im(\tau)>0\}$, 
and one defines 
$$\ord_\infty(f) = \min(n : a_f(n)\neq 0).
$$ 

The space of such forms $M_k$ is finite dimensional, spanned by the holomorphic Eisenstein series $E_k(\tau) = \frac 12 \sum_{\gcd(c,d)=1}(c\tau+d)^{-k}$  
($k\geq 4$) and the space of  cusp forms $S_k$, made up of forms which vanish at infinity ($\ord_\infty(f)\geq 1$). 
Writing 
\[
k=12 \ell +k', \quad k'\in \{0,4,6,8,10,14\}   
\]
then  $\ell = \dim S_k$  when $\ell\geq 1$,    and for any nonzero $f\in M_k$, we have $\ord_\infty(f)\leq \ell$.   

\begin{figure}[ht]
\begin{center}
\includegraphics[height=60mm]{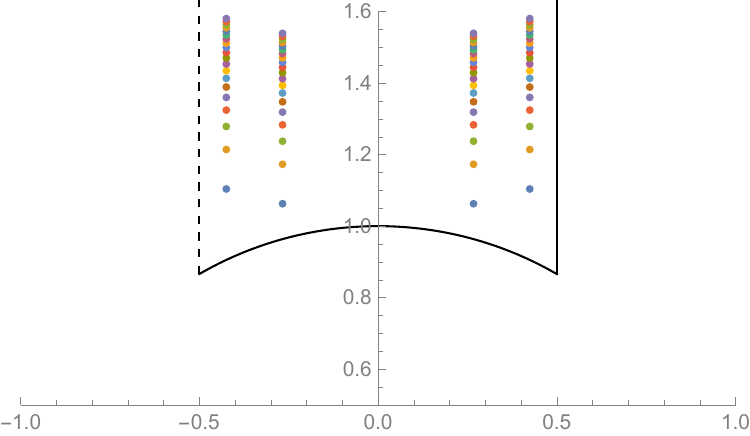}
\caption{The fundamental domain $\mathcal F$ and the points  $\frac{\sqrt{-1}}{2\pi} \log\left( 2k z_{4,r} \right)$, $r=1,\dots,4$ and $k= 1000 j$, $j=1,20$ where $z_{4,r}$ are the inverse zeros of $\Ex_4(t) = 1+t+t^2/2+t^3/6+t^4/24$.   }
\label{fig:level4}
\end{center}
\end{figure}

A nonzero modular form of weight $k$ has roughly  $k/12-\ord_\infty(f)$ zeros  (see \eqref{zeros formula}) in the fundamental domain  
\[
\mathcal F = \{\tau \in \HH: \re(\tau)\in [-1/2,1/2), |\tau|\geq 1\}
\]
where if $|\tau|=1$ then we take $\re(\tau)\leq 0$, see Figure~\ref{fig:level4}. 
The question we address is how are these zeros distributed as $k\to \infty$.

There are two main types of known results on zeros of modular forms:  
The first  of these originates in 1970,  when F.K.C. Rankin and P. Swinnerton-Dyer \cite{RanSD}  showed, by a remarkably simple argument,  that all the zeros of the Eisenstein series $E_k$  
lie on the bounding arc $\mathcal C=\{e^{it}:\pi/2\leq t\leq 2\pi/3\}$, and as $k\to \infty$ they become uniformly distributed there. 
Several authors have used the argument of Rankin and Swinnerton-Dyer as an ingredient in proving analogous results for other distinguished forms,  for instance 
  Duke and Jenkins \cite{DukeJenkins} studied certain ``gap forms'', see below.  
The second type of result concerns  cuspidal Hecke eigenforms. In this case,  the zeros are equidistributed in the fundamental domain with respect to the hyperbolic measure  
\cite{Rudnickzeros, HS}, see also \cite{GS, LMR} for further results on this case. 

Our goal is to present a different distribution result, for zeros of cusp form with a high order of vanishing at infinity. 
We investigate the  zeros of cusp forms with very high order of vanishing at infinity, by which we mean that we fix $D\geq 1$ and a bounding parameter $C>0$, and consider cusp forms with $\ord_\infty(f)=\ell-D$, 
\begin{equation}\label{def of special f}
f=q^{\ell-D}\left( 1+ y_f(1)q+\dots+y_f(D)q^D \right) + O\left( q^{\ell+1} \right)
\end{equation}
with coefficients $|y_f(j)|\leq C$, and take $k\gg 1$. 
A principal example are the cofinal elements of the ``Miller basis'' of $M_k$,   which are the unique elements $f_{k,m}$ of $M_k$ with $q$-expansion 
 \[
 f_{k,m} = q^m + O(q^{\ell+1}), \quad m=0,\dots, \ell .
 \]
 If we fix $D$, then $f_{k,\ell-D}=q^{\ell-D} + O(q^{\ell+1})$ is of the required form.

We show that for these forms, the zeros do not lie on the arc $\mathcal C$ as is the case for the Eisenstein series \cite{RanSD}, or for the ``gap form'' $f_{k,0}$ \cite{DukeJenkins},  
nor are the zeros equidistributed as is the case for cuspidal Hecke eigenforms of large weight. Instead, we find that the zeros asymptotically  lie on $D$ lines at height $\log 2k$. Precisely, let $\Ex_D(t)$ be the truncated exponential polynomial of degree $D$: 
\[
\Ex_D(t) = 1+t+\dots +\frac{t^D}{D!} 
\]
and denote by  $\{z_{D,r}\}$   the inverse zeros: $\Ex_D(t) = \prod_{r=1}^D\left(1-z_{D,r}t \right)$.
\begin{thm}\label{thm:zeros of f}
Fix $D\geq 1$, $C>0$, and let $f\in S_k$ be as in \eqref{def of special f}. 
Then the zeros $\tau_1,\dots, \tau_D$ of $f$ in the fundamental domain, suitably labeled, satisfy 
\[
\tau_r = \frac{\sqrt{-1}}{2\pi} \log\left( 2k z_{D,r} \right)  + O\left( \frac 1k \right) .
\]
\end{thm}
So in particular, the zeros of $f$ cluster around  the vertical lines 
$$\mathcal L_r=\{\re(\tau )  = -\frac{\arg(z_{D,r})}{2\pi} \} , \quad r=1,\dots, D$$  
where the argument is chosen so that $\arg(z_{D,r}) \in [-\pi,\pi)$, see Figure~\ref{fig:level4}. 


Our argument uses   Faber polynomials: Given a nonzero $f\in M_k$, the associated Faber polynomial $F_f(t)\in \C[t]$ is a polynomial of degree $D=\ell -\ord_\infty(f)$,  
uniquely determined by 
\[
\frac{f}{\Delta^\ell E_{k'}} = F_f(j)
\] 
where 
\[
j=\frac 1q + 744 +196884q+\dots 
\]
is Klein's absolute  invariant. The zeros  of $F_f$ are at  $j(\tau_r)$ where $\{\tau_r\}$ are the zeros of  $f/E_{k'}$. 
In \S~\ref{sec:Proof of Theorem thm:Faber} we show:
   \begin{thm}\label{thm:Faber}
 Let $f \in S_k$ be as in \eqref{def of special f} and $F_f(t)$ its Faber polynomial. Then 
 \[
\frac 1{(2k)^D} F_f(2k  t)  = \sum_{s=0}^D \frac 1{s!} \left( 1+O\left( \frac 1k\right) \right) t^{D-s}.
  \]
 \end{thm}
Noting that $ \sum_{s=0}^D \frac {t^{D-s}}{s!} = t^D\Ex_D(\frac 1t)$,  
we will then obtain the limit distribution of the zeros of $F_f(t)$:
 \begin{cor}\label{cor:zeros of F}
 As $k\to \infty$, the zeros $t_{1},\dots, t_{D}$ of $F_f(t)$ satisfy 
  \[
 t_{r} = 2k  \cdot z_{D,r} +O(1), \quad r=1,\dots,D,
 \]
 where $z_{D,1},\dots,z_{D,D}$ are the inverse zeros of $\Ex_D$. 
 \end{cor}
 In \S~\ref{sec:back to zeros of f} we will deduce    Corollary~\ref{cor:zeros of F} and Theorem~\ref{thm:zeros of f} from 
 Theorem~\ref{thm:Faber}.
 
 \bigskip
 
\noindent{\bf  Acknowledgements:} We thank Mikhail Sodin for a helpful discussion. This research was supported by the European Research Council (ERC) under the European Union's Horizon 2020 research and innovation programme (grant agreement No. 786758) and by the Israel Science Foundation (grant No. 1881/20). 
 
\section{Background on modular forms}

\subsection{Basic definitions} 
For an even integer $k\geq 0$, the space of modular forms $M_k$ consists of holomorphic functions on the upper half-plane $\HH=\{\tau=x+iy:y>0\}$ which transform under M\"{o}bius transformations from $\SL(2,\Z)$ as $f(\frac{a\tau+b}{c\tau+d}) = (c\tau+d)^k f(\tau)$ for all $\big(\begin{smallmatrix} a&b\\c&d \end{smallmatrix}\big)\in \SL(2,\Z)$, 
and are bounded as $\im(\tau)\to +\infty$. The subspace $S_k\subset M_k$ of cusp forms consists of those forms which vanish as 
$\im(\tau)\to +\infty$. A modular form has an expansion in terms of the nome $q=e^{2\pi i \tau}$:
\[
f(\tau) = \sum_{n= 0}^\infty a(n)q^n
\]
and is a cusp form when $a(0)=0$. 

Some examples are the  (normalized) Eisenstein series   
\[
E_k(z) = \frac 12 \sum_{\substack{ (c,d)\in \Z^2\\ \gcd(c,d)=1}} (cz+d)^{-k}
\] 
which have Fourier expansion 
\[
E_k(\tau) = 1 -\gamma(k) \sum_{n=1}^\infty \sigma_{k-1}(n)q^n
\]
where $\sigma_s(n) = \sum_{d\mid n} d^s$ are divisor sums, $\gamma(k)=2k/B_{k}$ with $B_k$  the Bernoulli numbers, see 
Table~\ref{table: gammak}.
 %
For instance, 
\[
E_4 = 1+240\sum_{n\geq 1} \sigma_3(n)q^n, \qquad E_6 = 1- 504\sum_{n\geq 1} \sigma_5(n)q^n.
\]
 An example of a cusp form is the modular discriminant $\Delta \in S_{12}$, the unique (up to multiple) cusp form of weight $12$, with Fourier expansion
\[
\Delta(\tau) = q\prod_{n=1}^\infty (1-q^n)^{24}  = q-24 q^2+252 q^3 \dots + . 
\]

\begin{table}[h]
\begin{center}\setlength{\tabcolsep}{10pt}   
\begin{tabular}{ |c|c|c|c|c|c|c|   }
\hline
 $k$ & $4$ & $6$ & $8$  & $10$ & $12$ & $14$   \\ [0.5ex] 
 \hline
 $\gamma(k)$ & -240  & $504$ & -480 & 264 & -65520/691 & 24  
 \\
 \hline
\end{tabular}
\end{center}
\caption{The numbers $  \gamma(k)=   {2k }/{B_{k}} $.}
\label{table: gammak}
\end{table}
 
We will also need  Klein's absolute invariant, the $j$-function
\[
j=   \frac{E_4^3}{\Delta }
\]
which is a modular function (weight zero), meromorphic at infinity, with $q$-expansion
\begin{equation}\label{q-expansion of j}
j 
=\frac 1q + 744 + 196884 q + 21493760 q^2 + \dots \in \Z[[q]] .
\end{equation}
The $j$-function gives an isomorphism $j:\SL(2,\Z)\backslash \HH \to \C$.  Any meromorphic modular form which is entire in the finite half-plane (its only possible poles are at infinity) is a polynomial in $j$. 

For $\tau$ in the standard  fundamental domain $\mathcal F$, $j(\tau)$ is real  if and only if $\tau$ lies on the boundary of $\mathcal F$  or on the imaginary axis, more precisely, $j$ maps the arc $e^{it}:t\in [\pi/2,2\pi/3]$ onto $[0,1728]$, the imaginary axis $\{iy:y\geq 1\}$ to $[1728,\infty)$ and the left boundary segment $\{-\frac 12 +iy:y> \sqrt{3}/2\}$ to the negative reals.

\subsection{Zeros}

Let $f\in M_k$  be a nonzero modular form of weight $k>0$. Then   the valence formula is: 
\begin{equation}\label{zeros formula}
\ord_\infty f + \sum_{\substack {z\in \mathcal F\\  z\not\sim i,\rho}} \ord_z f+\frac 12 \ord_{i}f + \frac 13 \ord_{\rho}f
=   \frac{k}{12}
\end{equation}
where $\rho=   (-1+\sqrt{-3})/2$ and the sum is over the  zeros   of $f$ in the fundamental domain $\mathcal F$ other than $\rho$ and $i$.  

From the valence formula \eqref{zeros formula} we see that if  $ k'\in \{ 4,6, 10 \}$ then the zeros of $E_{k'}$ are all simple; if $k'=8$ then there is a double zero at $\rho$ and no others; and  if $k'=14$ then   there is a simple zero at $i$ and a double zero at $\rho$.  We also see that, writing as before $k=12\ell+k'$, and form of weight $k$ has at least the same zeros as $E_{k'}$, so is divisible by $E_{k'}$ (for $k'=0$ we set $E_0=1$). 
We will use the term ``trivial zeros'' for the zeros of $f\in M_k$ arising from these symmetries.  The non-trivial zeros of $f$ are thus the zeros of the quotient $f/E_{k'}$.

\section{Faber polynomials}

\subsection{Definition}
To any modular form $f\in M_k$ we can associate a polynomial $F_f\in \C[t]$ so that 
\[
 f =\Delta^\ell E_{k'}  \cdot F_f(j) 
\]
(recall $k=12 \ell +k'$, $ k'\in \{0,4,6,8,10,14\}$). 
Indeed, the quotient $ f /\Delta^\ell E_{k'} $ is a meromorphic modular form, whose only possible poles are at infinity, 
hence  $ f /\Delta^\ell E_{k'} $   must be a polynomial in $j$, of  degree 
\[
D=\deg F_f = \ell - \ord_\infty(f) .
\]

By definition, for $ k'\in \{0,4,6,8,10, 14 \}$ we have $F_{E_{k'}}(t)=1$ and likewise for $f=\Delta^\ell$. Also by definition, multiplying $f$ by a power of $\Delta$ does not change the Faber polynomial: $F_{\Delta^mf} = F_f$.

The polynomial $F_f$ accounts for all the ``nontrivial'' zeros of $f$ (that is except for the common zeros with $E_{k'}$), in the sense that for these zeros $\tau$, we have $f(\tau)=0$ iff $F_f( j(\tau)) =0$. 



\subsection{Computation}
To compute the Faber polynomial $F_f(t)$, expand $ f /\Delta^\ell E_{k'} $ as  Laurent series in $q$ and then match the principal part with that of of a polynomial of degree $D$ in $j$:  
With $m:=\ord_\infty(f)$, expand
\begin{equation*}
\begin{split}
\frac {q^m}{\Delta^\ell E_{k'}} & =    \frac {q^m }{q^\ell\prod_{n=1}^\infty (1-q^n)^{24\ell}
 \left\{1-\gamma(k') \sum_{n\geq 1} \sigma_{k'-1}(n)q^n \right\}} 
 \\
 \\
 &= 
 \frac {A_k(0)+A_k(1)q +\dots +A_k(D)q^D}{q^D} +O(q)  .
\end{split}
\end{equation*}
with $A_k(0)=1$, and 
\[
f = q^m\sum_{n=0}^\infty y_f(n)q^n,  
\]
with $y_f(0)\neq 0$,  so that
\begin{equation}\label{expand f/DeltaE}
\begin{split}
\frac {f}{\Delta^\ell E_{k'}}  &=    
   \frac {A_k(0)+A_k(1)q +\dots +A_k(D)q^D}{q^D} \sum_{n=0}^\infty y_f(n)q^n   +O(q)  
   \\
  &= \sum_{s=0}^D q^{-s} \sum_{n=0}^{D-s} A_k(D-s-n)y_f(n) +O(q). 
\end{split}
\end{equation}

Further, expand
\[
j^r = \frac 1{q^r} + \frac{744 r}{q^{r-1}}+\dots = \sum_{s=0}^r c_{r,s}q^{-s} + O(q)
\] 
with $0\leq c_{r,s}\in \Z$, 
$$c_{r,r}=1, \quad c_{r,r-1} = 744 \cdot r $$
so that 
\begin{equation}\label{expand F(j)}
F_f(j) = \sum_{r=0}^D x_{D-r} j^r = \sum_{s=0}^Dq^{-s} \sum_{r=s}^D x_{D-r} c_{r,s} +O(q) .
\end{equation}

Comparing \eqref{expand f/DeltaE} with \eqref{expand F(j)} gives a system of equations for the coefficients $x_0,x_1,\dots x_D$ of the Faber polynomial:
\begin{equation}\label{eqs for x_r}
\sum_{r=s}^D  c_{r,s} x_{D-r}  = \sum_{n=0}^{D-s} A_k(D-s-n)y_f(n) , \quad s=0,1,\dots , D .
\end{equation}
For $s=D,D-1,D-2$    these are   
\[ 
 \begin{array}{lololololol}  
 x_0 &   &   &   &  & = &  y_f(0)     \\
744D      x_0  & + & x_{1} &   &   & = &  A_k(1)  y_f(0) &+ &y_f(1) &    \\
c_{D,D-2}  x_0   & + & c_{D-1,D-2}  x_{ 1}   & + &  x_{ 2}  & = &  A_k(2)   y_f(0)& +& A_k(1) y_f(1) &+&y_f(2) .
\end{array}   
\]


\subsection{Examples}
%
 We determine in this way the Faber polynomials $F_{k,m}(t)$ for some examples in the  Miller basis $f_{k,m} = q^m+O(q^{\ell+1})$. 
For instance, for the ``gap form'' $f_{k,0} = 1+O(q^{\ell+1})$, 
 \[
F_{24,0}(t) =125280 - 1440\; t + t^2   
\]
with zeros $93.0072$, $1346.99$, and 
\[
F_{36,0} = -27302400 + 965520 \; t - 2160 \; t^2 + t^3 
\]
with zeros $ 30.3029$, $582.232$, $1547.46 $. 
Duke and Jenkins \cite{DukeJenkins}  show the gap form $f_{k,0}$ has all it zeros on the arc 
$\{e^{i\theta}: \frac \pi 2\leq \theta \leq \frac{2\pi}{3}\}$, equivalently that the Faber polynomial has all its zeros in $[0,1728]$.

Here are some examples of $F_{k,\ell-D}$ when $ k=12\ell$ and $D=\ell-m$ is small: 
\[
\begin{split}
F_{12\ell,\ell-1}(t) &= t+(2k  -744 ).
\\ 
%
F_{12\ell,\ell-2}(t) &= t^2 +24 (\ell-62 )t + 36 (8 \ell^2 - 495 \ell +4438 )
\\
&=t^2 +\left( 2k -1488 \right) t + \left(\frac {(2k)^2}{2}-1485k +159768 \right) .
\\
%
F_{12\ell, \ell-3}(t) &= t^3 + 24 (-93 + \ell) t^2 + 36 (29721 - 991 \ell + 8 \ell^2) t 
\\ &\quad + 32 (-1152093 + 118990 \ell - 6669 \ell^2 + 72 \ell^3)
\\
&=t^3 + \left(2k-2232 \right)t^2 +\left( \frac{(2k)^2}{2} -2973 k  +1069956 \right)t   
\\ &\quad + \left( \frac{(2k)^3}{6} -1482 k^2  +\frac{951920}{3} k -36866976 \right) .
\end{split}
\]

 \section{ Proof of Theorem~\ref{thm:Faber}}\label{sec:Proof of Theorem thm:Faber}

 %

 
 Recall that we fix $D = \ell-m$, fix $C>0$, and consider   cusp forms of the shape 
 \[
 f=q^m\left(1+y_f(1)q+\dots +y_f(D) q^D \right) +O(q^{\ell+1})
 \]
 with bounded coefficients: 
 \[
  |y_f(1)|, \dots |y_f(D)| \leq C .
 \]
We will prove Theorem~\ref{thm:Faber}, which states that 
 \[
 \frac 1{(2k)^D} F_f(2k  t)  = \sum_{s=0}^D \frac 1{s!} \left( 1+O\left( \frac 1k\right) \right) t^{D-s}.
 \]
\begin{proof}
 Using 
 \[
 \frac 1{(1-x)^{N}}=  1+Nx+\dots = \sum_{r=0}^\infty \binom{N-1+r}{r}x^r
 \]
we obtain
 \[
 \begin{split}
 \frac 1{\prod_{n=1}^D(1-q^n)^{24\ell }} &= \prod_{n=1}^D \sum_{r_n=0}^D \binom{24\ell -1+r_n}{r_n} q^{n r_n} + O(q^{D+1})
 \\
 &=\sum_{r=0}^D q^r \sum_{r_1+2r_2+\dots +Dr_D=r}  \prod_{n=1}^D \binom{24\ell-1+r_n}{r_n} + O(q^{D+1}) .
 \end{split}
 \]
 As $\ell\to \infty$, the coefficient $B_{24\ell}(r)$ of $q^r$ in the above expansion is dominated by the contribution of the $d$-tuple $(r_1,\dots,r_D) = (r,0,\dots,0)$:
 \[
\begin{split}  
 B_{24\ell}(r)& = \sum_{r_1+2r_2+\dots +Dr_D=r}  \prod_{n=1}^D \binom{24\ell -1+r_n}{r_n} 
\\
 &=
 \sum_{r_1+2r_2+\dots +Dr_D=r}   \frac{(24\ell )^{r_1+r_2+\dots+r_D}}{r_1!r_2!\dots r_D!} \left(1+O\left(\frac 1\ell\right)\right)
\\ 
 &= \frac{(24\ell )^r}{r!} \left(1+O\left(\frac 1\ell\right)\right)
=  \frac{(2k)^r}{r!} \left(1+O\left(\frac 1k\right)\right) .
 \end{split}
 \]
 
 We next expand  $f/(q^m E_{k'})$ up to $O(q^{D+1})$:
 \[
 \frac {f}{q^m E_{k'}} = \frac {1+y_f(1)q+\dots +y_f(D)q^D}{1-\gamma(k')\sum_{n=1}^D \sigma_{k'-1}(n)q^n}+O(q^{D+1}) = \sum_{s=0}^D \alpha_s q^s  
 +O\left(q^{D+1} \right) 
 \]
 with $\alpha_0=1$ and  $$\alpha_s = O(1)   $$ 
 because we assume that the coefficients $y_f(1),\dots, y_f(D)$ are uniformly bounded as $k\to \infty$.
 
Multiplying by $\prod_{n=1}^D(1-q^n)^{-24\ell }$ gives that the Taylor polynomial of degree $D$  in the expansion of 
 $1/( E_{k'}\prod_{n=1}^D(1-q^n)^{24\ell}) $ is 
 \[
 \begin{split}
   \sum_{r=0}^D \frac{(24\ell )^r}{r!}q^r  \left( 1+O\left( \frac 1k \right) \right) \sum_{s=0}^D \alpha_s q^s 
& = \sum_{i=0}^D q^i \sum_{\substack{r+s=i\\r,s\geq 0}} \frac{(24\ell )^r}{r!} \alpha_s \left( 1+O\left( \frac 1k \right) \right)  
\\
&= \sum_{i=0}^D \frac{(24\ell )^i}{i!} \left(1+O\left( \frac 1k \right) \right)q^i  
\end{split} 
\]
since $\alpha_0=1$ and $\alpha_s = O(1)$. 
 Therefore, the terms up to $O(q)$ of $f /(\Delta^\ell E_{k'})$ are 
 \[
  \frac 1{q^D}\sum_{r=0}^D \frac{(24\ell)^r}{r!} \left( 1+O\left( \frac 1k \right) \right) q^r    .
  \]
 Finally, replacing $24\ell $ by $2k=24\ell +O(1)$ we obtain
  \[
\frac{f }{\Delta^\ell E_k'} 
= \frac 1{q^D}\sum_{r=0}^D A_{k}(r)q^r + O(q)  
\]
with
\[
 A_{k}(r) = \frac{(2k)^r}{r!} \left( 1+O\left( \frac 1k \right) \right), \quad  r=0,\dots ,D .
\]

 Now compare with $F_f(j)=j^D+x_1j^{D-1}+\dots +x_D$: Using the equations \eqref{eqs for x_r} gives
 \[
 744 D + x_1 = 2k \left( 1+O\left( \frac 1k \right) \right)
  \]
 which says that 
 \[
  x_1 = 2k  +O(1) .
 \] 
Next, we have 
 \[
 c_{D,D-2} + c_{D-1,D-2} x_1 + x_2 = \frac{(2k)^2}{2!} + O(k)
 \]
 and since $x_1 = O(k)$ we get
 \[
 x_2 =  \frac{(2k)^2}{2!} + O(k) .
 \]
 Continuing, we assume by induction that 
 \[
 x_1 = O(k), \dots, x_{s-1} = O(k^{s-1})
 \]
 and then obtain
 \[
 c_{D,D-s} + \sum_{i=1}^{s-1} c_{D-i,D-s} x_i + x_s =  \frac{(2k)^s}{s!} + O(k^{s-1})
 \]
which gives
\[
 x_s =  \frac{(2k)^s}{s!} + O(k^{s-1}) .
\]

Thus
\[
\begin{split}
\frac 1{(2k)^D}F_f(2k t)  &= \frac 1{(2k)^D}\sum_{s=0}^D  \frac{(2k)^s}{s!}(2k t)^{D-s} \left( 1+O\left( \frac 1k \right) \right) 
\\
&=   \sum_{s=0}^D \frac{1}{s!}  \left( 1+O\left( \frac 1k \right) \right) t^{D-s} .
\end{split}
\]
as claimed. 
\end{proof}

\section{Back to zeros of modular forms}\label{sec:back to zeros of f}
Having at hand the convergence of the coefficients of the renormalized Faber polynomials $F_f(2kt)/(2k)^D$ to those of $ t^D \Ex_D(1/t)$, we can deduce convergence of zeros.

\subsection{Proof of Corollary~\ref{cor:zeros of F}}
We set
\[
g_k(z) =\frac{ F_f(2k \cdot z)}{(2k)^D} , \quad g(z) = z^d\Ex(\frac 1z) = \sum_{r=0}^D \frac{z^{D-r}}{r!} . 
\]
The zeros of  $g(z)$ are simple, as follows from the corresponding fact for $\Ex_D(z)$, which is in fact irreducible; see 
\cite{Zemanyan} for a survey.

Moreover, 
\[
\begin{split}
g_k(z) &= \sum_{r=0}^D \frac{(2k)^r}{r!} \left( 1+O\left( \frac 1k \right) \right) 
\frac{ (2k z)^{D-r}}{(2k)^D} 
\\ 
&=   \sum_{r=0}^D\frac{t^{D-r}}{r!}\left( 1+O\left(\frac 1k \right) \right) 
\end{split}
\]
Therefore, we deduce that for $k\gg 1$, the zeros $z_{k,1},\dots, z_{k,D}$ of $g_k$ are simple and converge to the zeros $z_1,\dots, z_D$ of $g$, with a rate
\begin{equation}\label{eq:Approximation of zeros}
z_{k,r} = z_r + O(\frac 1k). 
\end{equation}
This follows for instance from \cite[Appendix A, Theorem]{Ostrowksi} which states that given monic polynomials 
$f(t)=t^D+\sum_{\nu=1}^D a_\nu t^{D-\nu}$ with zeros $x_1,\dots, x_D$ and $g(t) =t^D+\sum_{\nu=1}^D b_\nu t^{D-\nu}$ with zeros $y_1,\dots, y_D$,  then possibly after relabeling the zeros,  we have a bound on their differences
\[
\max_{\nu=1,\dots,D} |x_\nu-y_\nu | \leq 2D \left( \sum_{\nu=1}^D |a_\nu-b_\nu| \;\Gamma^{D-\nu} \right)^{1/D}
\]
 where 
 \[
 \Gamma = \max_{\nu=1,\dots,D}( |a_\nu|^{1/\nu}, |b_\nu|^{1/\nu} ) .
\]
In our case, taking $f=g_k$ for $k\gg 1$, and $g(z) = \sum_{\nu=0}^D z^{D-\nu}/\nu!$, we clearly have $|a_\nu-b_\nu|\ll 1/k$ and 
$\Gamma =O(1)$  and so we obtain \eqref{eq:Approximation of zeros}

The zeros of $F_f$ are $t_{k,r} =2k z_{k,r}$, hence \eqref{eq:Approximation of zeros} implies that they satisfy
\[
t_{k,r} = 2k z_r =2k z_r+O(1)
\] 
which proves   Corollary~\ref{cor:zeros of F}. \qed

\subsection{Proof of Theorem~\ref{thm:zeros of f}}
The nontrivial zeros $\tau_1,\dots,\tau_D$ of $f$, are the zeros of $F_f\left(j\left(\tau\right)\right)$, so their $j$-values satisfy 
\[
 j\left(\tau_r\right)  = 2k z_{D,r} + O(1), \quad r=1,\dots, D.
\] 
Therefore these $j$-values tend to infinity, and in terms of the nome $q_r=e^{2\pi i \tau_r}$, 
\[
\frac 1{q_r} +O(1) = j(\tau_r) = 2k z_{D,r}  +O(1) .
\]
Hence 
\[
e^{ -2\pi i \tau_r} =  2k z_{D,r}  +O(1) 
\]
giving
\[
\tau_r = \frac i{  2\pi i} \log \left( 2k z_{D,r}  \right) + O\left(\frac 1k \right)
\]
as claimed in Theorem~\ref{thm:zeros of f}. \qed

\end{document}